\documentclass[12pt]{amsart}
\usepackage{amssymb}
\usepackage{amsfonts}
\usepackage{amsmath}
\usepackage{array}

\usepackage{latexsym}

\newtheorem{example}{Example}[section]
\newtheorem{note}[example]{Note}
\newtheorem{theorem}[example]{Theorem}

\newtheorem{corollary}[example]{Corollary}

\newtheorem{definition}[example]{Definition}
\newtheorem{proposition}[example]{Proposition}

\newtheorem{lemma}[example]{Lemma}

\def\Proof{\noindent \it Proof -- \rm}
\def\qed{\hspace{3.5mm} \hfill \vbox{\hrule height 3pt depth 2 pt width 2mm}
\bigskip}

\def\<{\langle}
\def\>{\rangle}

\def\C{\operatorname{\mathbb C}}

\def\N{\operatorname{\mathbb N}}

\def\G{{\bf G}}

\def\SG{{\mathfrak S}}
\def\H{{\mathcal H}}
\def\t{{\bf t}}

\def\Sym{{\bf Sym}}

\def\SymN{\hbox{\bf Sym}(N)}
\def\vtr#1{\vrule height 0mm depth #1mm width 0mm}
\def\wtr#1#2{\vrule height #2mm depth #1mm width 0mm}

\def\Ff{{\mathfrak F}}
\def\G{{\mathfrak G}}

\def\Tabvrule{\vrule width-0.4pt}       
\def\Tabhrule{\hrule \hrule height-0.4pt} 
\def\Tabstrut{\vrule height2.2ex 
                     depth0.8ex  
                     width0ex    
\relax}

\def\PasCase#1{\omit%
            $\vcenter{\hbox {\vbox to 0.4pt{}}
               \hbox{\makebox[3ex]{\Tabstrut$#1$}}}%
               \Tabvrule$}
\def\PasCasePoint{\PasCase{\cdot}}
\def\DessinCarre#1{%
    \vcenter{\hbox{}\hrule
             \hbox{\vrule\makebox[3ex]{\Tabstrut$#1$}\vrule}\Tabhrule}%
             \Tabvrule}
\def\GenRuban#1{\vcenter{\halign{&$\DessinCarre{##}$\cr#1}}\egroup}

\def\sTabvrule{\vrule width-0.4pt}
\def\sTabhrule{\hrule \hrule height-0.4pt}
\def\sTabstrut{\vrule height1.6ex depth0.6ex width0ex \relax}
\def\sDessinCarre#1{%
    \vcenter{\hbox{}\hrule
             \hbox{\vrule\makebox[2.3ex]%
                  {\sTabstrut$\scriptstyle#1$}\vrule}\sTabhrule}%
             \sTabvrule}
\def\sGenRuban#1{\vcenter{\halign{&$\sDessinCarre{##}$\cr#1}}\egroup}

\def\ruban{%
  \bgroup
  \let\ =\omit
  \let\\=\cr
  \let\x=\times
  \let\.=\PasCasePoint
  \offinterlineskip
  \GenRuban}

\def\sruban{%
  \bgroup
  \let\ =\omit
  \let\x=\times
  \let\\=\cr
  \offinterlineskip
  \sGenRuban}

\title{Higher order peak algebras}

\author[D.~Krob and J.-Y.~Thibon]%
{Daniel Krob and Jean-Yves Thibon}

\address[D. Krob]{LIX (CNRS), Ecole Polytechnique,\\
Route de Saclay,\\ 91128 Palaiseau Cedex \\
FRANCE}
\address[J.-Y. Thibon] {Institut Gaspard Monge, Universit\'e de Marne-la-Vall\'ee \\
5 Boulevard Descartes \\Champs-sur-Marne \\77454 Marne-la-Vall\'ee cedex 2 \\
FRANCE}
\email[D.Krob]{dk@lix.polytechnique.fr}
\email[Jean-Yves Thibon]{jyt@univ-mlv.fr} 
\date{}

\begin{document}

\begin{abstract}
Using the theory of noncommutative symmetric functions, we introduce
the {\it higher order peak algebras}
$(\SymN)_{N\geq 1}$, a sequence of graded Hopf algebras
which contain
the  descent algebra and the usual peak algebra
as initial cases  ($N = 1$ and $N = 2$). We
compute their Hilbert series, introduce and study several combinatorial
bases, and establish various algebraic identities
related to the multisection of formal power series with noncommutative
coefficients.
\end{abstract}

\maketitle

\section{Introduction}

The peak algebra of the symmetric group $\SG_n$ is a subalgebra of its
descent algebra $\Sigma_n$, spanned by sums of permutations having the same
peak set, a certain subset of the descent set.
The direct sum of these peak algebras turns out to be a Hopf subalgebra
of the direct sum of all descent algebras, which can itself
be identified with
$\Sym$, the Hopf algebra of noncommutative symmetric functions.

The peak Hopf algebra $\Pi$ has been introduced (somewhat implicitely, and independently)
in two papers published in 1997. In \cite{St}, Stembridge introduced the dual Hopf
algebra as a subalgebra of quasi-symmetric functions, while in
\cite[Sec.~5.6.4]{NCSF2}, Prop.~5.41 implied the existence of $\Pi$ and many
of its properties (see \cite{BHT}). Since then, a larger peak algebra,
related to descent algebras of type $B$ and $D$, has been discovered \cite{ABN}.

Actually, the results of \cite{NCSF2} implied more than
the classical peak algebra. The peak algebra was
obtained as the image of $\Sym$ under a certain endomorphism $\theta_q$
for $q=-1$. For generic values of $q$, $\theta_q$ is an automorphism, but
it degenerates at (non trivial) roots of unity. In this way, we obtain
an infinite family of Hopf subalgebras of $\Sym$, denoted by $\SymN=\theta_q(\Sym)$
when $q$ is a primitive $N$th root of unity. Then $\Sym$ and $\Pi$ correspond
to $N=1$ and $N=2$, respectively.

In the following, we shall investigate the Hopf algebras $\SymN$ for arbitrary
$N$. Our main results consist in the introduction of several linear bases,
in which the various natural operations of $\SymN$ admit interesting combinatorial
expressions. These bases allow us to obtain algebraic identities
generalizing those involving the noncommutative tangent in \cite{NCSF1,BHT}.

Finally, let us mention that the
commutative images (in ordinary symmetric functions)
of the first two algebras are known to be the Grothendieck rings of Hecke algebras
$H_n(q)$ at $q=1$ (symmetric groups, same as generic case) and $q=-1$.
In general, the commutative image of the $N$-th peak algebra is the
Grothendieck ring of the tower $H_n(\zeta)$, where $\zeta$ is
a primitive $N$-th root of unity (symmetric functions not
depending on power sums $p_m$ with $N|m$, cf. \cite{LLT}). 
It would be interesting to find
an interpretation of the higer order peak algebras in this context.



\section{Preliminaries}

Our notations will be those of \cite{NCSF1,NCSF2,NCSF3,NCSF4,NCSF5,NCSF6}. We
recall here the most important ones in order to make the paper self-contained.


\subsection{Compositions and permutations}
\label{subsec:comp_perm}

A {\it composition} of $n$ is a sequence $I =(i_1,i_2,\cdots,i_r)$ 
of positive integers summing to
$n$. The {\it length} of $I$ is  $\ell(I) = r$.
A composition $I$ can be represented by
a skew Ferrers diagram called a {\it ribbon diagram} of shape $I$, e.g.,

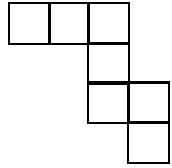
\setlength{\unitlength}{0.3pt}
\centerline{
\begin{picture}(250,255)(0,-165)
\put(0,0){\framebox(50,50){}}
\put(50,0){\framebox(50,50){}}
\put(100,0){\framebox(50,50){}}
\put(100,-50){\framebox(50,50){}}
\put(100,-100){\framebox(50,50){}}
\put(150,-100){\framebox(50,50){}}
\put(150,-150){\framebox(50,50){}}
\end{picture}
}
\noindent is a ribbon diagram of shape $I = (3,1,2,1)$.
The {\it conjugate composition} $I^{\sim}$ is 
associated with the conjugate ribbon diagram
in the sense of skew Ferrers diagrams. On our example,
$I^{\sim} = (2,3,1,1)$.

\smallskip
One associates with a composition $I = (i_1,i_2,\cdots,i_r)$ of
$n$ the subset $D(I)$ of $[1,n\!-\!1]$  defined by
$D(I) = \{\,i_1,i_1+i_2,\, \cdots\,,i_1+i_2+\cdots+i_{r-1}\,\}$.
Compositions of  $n$ 
are ordered by {\it reverse refinement}, denoted by $\preceq$,
and defined by  $I \preceq J$ iff $D(I) \subseteq D(J)$.

\smallskip
One associates with a permutation $\sigma$ of $\SG_n$ two compositions
$D(\sigma)$ and $P(\sigma)$ respectively called the {\it descent composition}
and the {\it peak composition} of $\sigma$. The descent composition $D(\sigma)$
of the permutation $\sigma$ encodes its descent set
${\mathcal D}(\sigma) = \{\, i \in [1,n\!-\!1], \, \sigma(i) > \sigma(i\!+\!1)\,\}$,
i.e., is characterized by $D(D(\sigma)) = {\mathcal D}(\sigma)$.
The peak composition $P(\sigma)$ of $\sigma$ encodes its peak
set
${\mathcal P}(\sigma) = \{\, i \in [2,n\!-\!1], \, \sigma(i\!-\!1) < \sigma(i) > \sigma(i\!+\!1)\,\}$,
that is,  $D(P(\sigma)) = {\mathcal P}(\sigma)$.
For example $D(265341) = (2,1,2,1)$, ${\mathcal D}(265341) = \{\, 2,3,5\,\}$,
$P(265341) = (2,3,1)$ and ${\mathcal P}(265341) = \{\, 2,5\,\}$.

\smallskip
While each subset of $[1,n\!-\!1]$ is obviously the descent set of
some permutation $\sigma \in \SG_n$, a subset $P$ of $[1,n\!-\!1]$
is the peak set of a permutation $\sigma \in \SG_n$ if and only
if 
\begin{equation}
  \left\{ \
  \begin{array}{l}
    1 \notin P \, , \wtr{3}{4} \\
    \forall \, i \in [2,n\!-\!1], \, i \in P \,\Rightarrow \, i-1 \notin P\ .
    \vtr{2.5}
  \end{array}
  \right.
\end{equation}
The corresponding compositions are  called 
{\it peak compositions} of  $n$.
A peak composition of  $n$ is 
one whose only entry allowed to be equal to $1$ is the last
one.  It is easy to check that
the number of  peak compositions of  $n$ is
$f_n$, where $f_n$ is the $n$-th Fibonacci
number, with the convention $f_0 = f_1 = f_2 = 1$
and $f_{n+2} = f_{n+1} + f_n$.

\smallskip
Finally, the peak set $P(I)$ of a 
composition $I = (i_1,i_2,\cdots,i_r)$ 
is defined as follows. If $D(I) = \{d_1,d_2,\cdots,d_{r-1}\}$ with
$d_1 \leq \cdots \leq d_{r-1}$, then $P(I)$ is obtained from
$D(I)$ by removing each $d_i$ such that $d_i - d_{i-1} = 1$ 
(with the convention $d_0 = 0$). For example, if
$I = (1,3,1,2)$,  $D(I) = \{1,4,5\}$ and  $P(I) = \{4\}$.


\subsection{Noncommutative symmetric functions}

The Hopf algebra of noncommutative symmetric functions, 
introduced in \cite{NCSF1}, is the free associative
$\C$-algebra $\Sym = \C\<S_1,S_2,\ldots,\>$ over an infinite
sequence of noncommuting indeterminates $(S_k)_{k\geq 1}$  called the
{\it noncommutative complete symmetric functions}. 
It is graded by ${\rm deg\,}S_n=n$, and endowed with the coproduct
\begin{equation}
\Delta(S_n) = \ \sum_{k=0}^n \ S_k \otimes S_{n-k} \ .
\end{equation}
Its graded dual is Gessel's algebra of quasi-symmetric functions.
Various applications are discussed in the series
\cite{NCSF1,NCSF2,NCSF3,NCSF4,NCSF5,NCSF6}.

\smallskip
Let $\sigma(t)$ be the generating series of  the 
$(S_k)_{k\geq 0}$, with the convention $S_0 = 1$, 
\begin{equation}
\sigma(t) = \ \sum_{k\geq\, 0} \ S_k\,t^k \ .  \wtr{0}{5.5}
\end{equation}
The {\it noncommutative power sums of the first kind}, denoted by $({\Psi_k})_{k\geq 1}$,
are  the coefficients of the  series 
\begin{equation}
\psi(t) = \ \sum_{k\ge 1}\ \Psi_k\, t^{k-1}
        = \sigma(t)^{-1}\ {d\over dt}\,\Big(\, \sigma(t)\,\Big) \ .  \wtr{0}{5.5}
\end{equation}
The {\it noncommutative ribbon Schur functions} $R_I$
are characterized by  the  equivalent relations
\begin{equation}
S^I =\ \sum_{J\preceq \,I}\ R_J\ ,
\quad
R_I = \ \sum_{J\preceq\, I}\ (-1)^{\ell(I)-\ell(J)}\, S^J \wtr{0}{5.5}
\end{equation}
where $S^I = S_{i_1}\, S_{i_2} \,\cdots\, S_{i_r}$.

\smallskip
The {\it commutative image} of a noncommutative symmetric function $F$ is the
ordinary symmetric function $\chi(F)$  obtained through the
algebra morphism  $\chi(S_n) = h_n$,
the  complete symmetric function
(notation as in \cite{McD}). The commutative image of
$\Psi_n$ is then the  power sum $p_n$. Similarly,
the commutative image of $R_I$ is the ordinary ribbon Schur function of shape
$I$. We recall that ribbon Schur functions were introduced by McMahon (see
\cite{MM}, t. 1, p. 200, where they are denoted by $h_I$).


\subsection{Relation with the descent algebra of the symmetric group}
\label{subsec:rel_des_algebra}

The sum in the group algebra $\C[\SG_n]$ of all permutations
$\sigma$ whose descent composition $D(\sigma)$ is equal to $I$ is
denoted by $D_{=I}$. The  $D_{=I}$ with
$|I|=n$ form a basis of a subalgebra of $\C[\SG_n]$, introduced by
Solomon in \cite{So}, which is denoted $\Sigma_n$ and called the
{\it descent algebra} of $\SG_n$. One can define  an isomorphism
$\alpha$ of graded vector spaces
\begin{equation}
\alpha : \, \Sym = \ \bigoplus_{n\ge\, 0} \ \Sym_n
\longrightarrow \
\Sigma = \ \bigoplus_{n\ge\, 0} \ \Sigma_n\
\wtr{0}{5.5}
\end{equation}
by  $\alpha(D_{=I }) = R_I$. Then,
$\alpha(S^I) = D_{\preceq \, I}$,  the sum in $\C[\SG_n]$
of all permutations $\sigma$ whose descent composition $D(\sigma)$ is
$\preceq I$. 

\smallskip
The direct sum $\Sigma$ can be turned into an algebra by extending the natural
products of its components $\Sigma_n$ by setting $x\,y = 0$ for $x \in \Sigma_p$
and $y \in \Sigma_q$ with $p \not= q$. We can can then define the {\it internal
product}, denoted by $*$, on $\Sym$ by requiring that $\alpha$  be an anti-isomorphism
of algebras. In other words, the internal product $*$ of $\Sym$ is defined by
\begin{equation}
F * G = \alpha^{-1}(\alpha(G) \, \alpha(F))\,.
\end{equation}
%


\subsection{The $(1-q)$-transform} \label{sec:1mq}
\label{subsec:the_1-q_transformation}

In the commutative case, the {\it $(1\!-\!q)$-transform} is
the endomorphism $\vartheta_q$ of the algebra $Sym$ of 
commutative symmetric functions  defined on 
power sums by  $\vartheta_q(p_n) =
(1-q^n)\,p_n$. 
Our terminology comes from
the $\lambda$-ring notation, which allows to write it as
$\vartheta_q(f(X)) = f((1\! -\! q)\,X)$ (cf.
\cite{Knu,LS}). One just has to point out here the (traditional) abuse of
notation in using the same minus sign for  scalars and in the
$\lambda$-ring, though these operations are completely
different.

\smallskip
In \cite{NCSF2}, this $\lambda$-ring formalism was extended to the algebra of
noncommutative symmetric functions. A consistent definition of $\vartheta_q(F)
= F((1\! -\! q)\,A)$ for $F \in \Sym$ was proposed and its fundamental properties
were  obtained. The first step in the construction of this noncommutative
version consists in defining the noncommutative complete symmetric functions
$(S_n((1\!-\!q)\,A)_{n\geq\, 0}$ by their generating series
\begin{equation}
\sigma_t((1\!-\!q)\,A)
=
\ \sum_{n\ge 0} \ S_n((1\!-\!q)\,A)\,t^n
=
\sigma_{-q\,t}(A)^{-1}\, \sigma_t(A) \, ,
\end{equation}
where $\sigma_t(A)$ is the generating series of the complete
functions $S_n(A)$ (usually denoted by $\sigma(t)$ when $A$
plays no role).
The noncommutative $(1\!-\!q)$-transform is the algebra endomorphism
$\vartheta_q$ of $\Sym$  defined by  $\theta_q(S_n) = S_n((1\!-\!q)A)$.
One can show that 
\begin{equation}
\vartheta_q(F(A)) = F((1\!-\!q)\,A) = F(A) * \sigma_1((1\!-\!q)\,A)
\end{equation}
for  $F(A) \in \Sym$, where $*$ is the internal product.

\smallskip
The most important property of $\vartheta_q$ obtained in \cite{NCSF2} is its
diagonalization. We recall first that a noncommutative symmetric function
$\pi \in \Sym$ is called a {\it Lie idempotent} of order $n$ if $\pi$ is an
homogeneous noncommutative symmetric function of $\Sym_n$ 
in the primitive Lie algebra $L(\Psi)$, generated by  $(\Psi_n)_{n\geq 1}$,
and which is idempotent for the
internal product of $\Sym$, or equivalently which has a commutative image
equal to $p_n/n$ (see \cite{NCSF1}). It is shown in \cite{NCSF2} that there
is a unique family $(\pi_n(q))_{n\geq 1}$ of Lie idempotents
(with $\pi_n(q) \in \Sym_n$) possessing the characteristic property
\begin{equation} \label{equ:charproppi}
\vartheta_q(\pi_n(q)) = (1-q^n)\,\pi_n(q)\,.
\end{equation}
These Lie idempotents were used in \cite{NCSF2} for
describing the structure of $\vartheta_q$, which is a semi-simple endomorphism of
$\Sym_n$. The eigenvalues of $\vartheta_q$ in $\Sym_n$ are 
\begin{equation}
p_\lambda(1-q) = \ \prod_{i=1}^r \ (1-q^{\lambda_i})
\end{equation}
where $\lambda = (\lambda_1,\dots,\lambda_r)$ runs over all partitions of
$n$. The projector on the eigenspace
${\mathcal E}_{\lambda}$  associated with the eigenvalue $p_\lambda(1-q)$
is the endomorphism
$F \longrightarrow F * E_\lambda(\pi)$
of $\Sym_n$, where $E_\lambda(\pi)$ is 
obtained through a construction introduced in \cite{NCSF2}, which we briefly
recall here. Consider first the  decomposition of $S_n$ in the
multiplicative basis of $\Sym$ 
constructed from the  $\pi_n(q)$, i.e.,
\begin{equation}
S_n = \ \sum_{|I|=n} \ c_I(q)\,\pi^I(q)
\end{equation}
where $\pi^I(q) = \pi_{i_1}(q) \dots \pi_{i_r}(q)$.
Then, the projector $E_\lambda(\pi)$ on ${\mathcal E}_\lambda$ is given by
\begin{equation}
E_\lambda(\pi) = \ \sum_{I^\downarrow = \lambda} \ c_I(q)\,\pi^I(q)
\end{equation}
where $I^\downarrow$ denotes the  partition  obtained by reordering
the components of  $I$. It can be shown that the eigenspace
${\mathcal E}_{\lambda}$ is spanned by the $\pi^I(q)$ such that
$I^\downarrow = \lambda$ (see the proof of Theorem 5.14 of \cite{NCSF2}).

\smallskip
The following  formula can  be easily deduced from these results.

\begin{proposition} \label{prop:detthetaq}
Let $q \in \C$ be a complex number. Then we have :
\begin{equation} \label{equ:detthetatq}
{\det}\,\vartheta_q|{\Sym_n}
=
\left(\ \prod_{i=1}^{n-1} \ (1-q^i)^{(n-i+3)\,2^{n-i-2}} \, \right) \, (1-q^n) \ .
\end{equation}
\end{proposition}

\Proof Since the determinant of an endomorphism is the product of its eigenvalues,
we  get 
$$
{\det}\,\vartheta_q|{\Sym_n}
=
\ \prod_{{I = (i_1,\dots,i_r)  \wtr{1}{2} \atop |I| = n}}
\left( \ \prod_{k = 1}^r \ (1-q^{i_k}) \, \right)
=
\ \prod_{i=1}^n \ (1-q^i)^{\, c(n,i)}
$$
where $c(n,i)$ is  the number of times the integer $i$ appears as a component
of a composition of $n$. Formula (\ref{equ:detthetatq}) can therefore be established
if we can find a closed expression for this number. To this purpose, let
$$
F_i(x,y)
=
\ \sum_{k\ge 0} \ (x + x^2 + \dots + x^{i-1} + x^i\,y + x^{i+1} + x^{i+2} + \dots \,)^k
=
\ \sum_{k\ge 0} \ \left(\, {x \over 1-x} - x^i + x^i\,y \, \right)^k
\ .
$$
Then,
$$
\sum_{n\ge 0} \ c(n,i)\,x^n
=
\ {d \over dy}(F_i(x,y))\left.\right|_{\, y = 1} \,
=
x^i \, \left( \ \, \sum_{k\ge 0} \ k\, \Big(\, {x \over 1-x}\,\Big)^{k-1} \ \right)
=
x^i \, {1 - 2\,x + x^2 \over (1-2\,x)^2} \ .
$$
This shows  that  $c(n,i) = c(n\!-\!i\!+\!1,1)$, and
that $c(n,1)$ is equal to $1$ if $n = 1$ and to $(n+2)\,2^{n-3}$ for $n \geq 2$, 
whence the result.
\qed

\begin{note} \label{note:injection}
{\rm
If $q$ is  a complex number,  Proposition \ref{prop:detthetaq} shows that
$\vartheta_q$ is an automorphism of $\Sym$ if and only if $q$ is not a root of unity.
The higher order peak algebras $\SymN$ arise precisely
when we specialize $q$ to a primitive $N$-th
root of unity.
}
\end{note}


\subsection{The peak algebra}
\label{subsec:rel_peak_algebra}

The sum in the group algebra $\C[\SG_n]$ of all permutations
$\sigma$ whose peak composition is equal to $I$ is denoted by
$P_{=I}$. It follows from the results
of \cite{NCSF2} that the family  $(P_{=I})$, where $I$
runs over peak compositions of  $n$, form a basis
of a subalgebra of $\C[\SG_n]$, 
first mentioned explicitly by Nyman in \cite{Ny}.
It is denoted by $\Pi_n$ and called the \begin{it}peak
algebra\end{it} of $\SG_n$. Note that $\Pi_n$ is a subalgebra of
$\Sigma_n$. Indeed, for every peak composition $I$ of  $n$,
one has
$$
P_{=I} = \sum_{P(J)=D(I)}D_{=J} \ ,
$$
where $J$ runs over all compositions of $n$. Using the isomorphism
$\alpha$ defined in Section \ref{subsec:rel_des_algebra}, we can
therefore define the {\it peak functions} $\Pi_I = \alpha(P_{=I})$
in $\Sym_n$ by 
$$
\Pi_{I} = \sum_{P(J)=D(I)}R_J
$$
for all peak compositions $I$. The space spanned by all the peak
classes in $\Sym_n$ is called
$\mathcal{P}_n$. It is obviously isomorphic to $\Pi_n$, and so
$\mathcal{P} = \bigoplus_{n \geq 0} \mathcal{P}_n$ is
isomorphic to $\Pi = \bigoplus_{n \geq 0} \Pi_n$.

\smallskip
We will now explain the connection between these constructions
and the $(1-q)$-transform $\theta_q$ introduced in
Section \ref{subsec:the_1-q_transformation}
(see also \cite{BHT}).
Let us first recall the following result which expresses
$\theta_q(R_I)$ in the ribbon basis of $\Sym$ (cf. equation (121)
of \cite{NCSF2}),  
\begin{equation} \label{equ:dev_ri}
R_I((1-q)A) = \sum_{P(J) \subseteq D(I) \triangle (D(I) +
  1)}(1-q)^{hl(J)}(-q)^{b(I,J)}R_J(A)\,,
\end{equation}
Here $I$, $J$ are compositions of $n$, $hl(J) = |P(J)| +
1$ and $b(i,j)$ is an integer which will not bother us at this
point since we set $q = -1$ from now on. Rewriting this equation
in terms of peak classes of $\Sym_n$ gives 
$$
\theta_{-1}(R_I) = \sum_{D(J) \subseteq D(I) \triangle (D(I) + 1)}
2^{|D(J)| + 1}\Pi_{J} \ ,
$$
where $J$ runs over  peak compositions of  $n$. Hence,
the image of $\Sym_n$ by
$\theta_{-1}$, which we will denote by $\widetilde{\Sym_n}$, is
contained in the subspace $\mathcal{P}_n$ of $\Sym_n$ spanned by
the peak classes. In fact,
$$
\widetilde{\Sym_n} = \mathcal{P}_n\,,
$$
since
$\widetilde{\Sym_n}$ and $\mathcal{P}_n$ have the same dimension.
Indeed, recall that the dimension of $\Pi_n$ (and therefore of
$\mathcal{P}_n$) is equal to the Fibonacci number $f_n$ (as shown
in Section \ref{subsec:comp_perm}). Thanks to equation
(\ref{equ:charproppi}), we know that the elements
$$
\pi^I(-1) = \pi_{i_1}(-1)\hspace{0.1cm}\pi_{i_2}(-1)
\hspace{0.1cm} \cdots \hspace{0.1cm}\pi_{i_r}(-1)
$$
where $I=(i_1,\cdots,i_r)$ runs over all compositions of $n$ into
odd parts, form a basis of $\widetilde{\Sym_n}$, since
$\theta_{-1}(\pi_k(-1)) = 0$ if and only if $k$ is even. But the
number of compositions of $n$ with odd parts is easily seen to be $f_n$.
Therefore,
$$
\mathcal{P} = (\widetilde{\Sym} := \bigoplus_{n \geq 0}
\widetilde{\Sym_n})\ .
$$
%


\section{Higher order peak algebras}

\subsection{Definition and first properties}

Let $N \geq 1$ and let $\zeta \in \C$ be a primitive $N$-th root of
unity. We set $\Theta_q = \frac{\theta_q}{1-q}$, and we
denote  by $\Theta_{\zeta}$ the endomorphism of $\Sym$ 
defined by 
\begin{equation}
\Theta_{\zeta}(S_n) = \frac{S_n((1-q)A)}{1-q}\bigg{|}_{q = \zeta}
\end{equation}
\noindent
for  $n \geq 1$. Due to the fact that $\Theta_{\zeta}$ is not
invertible (cf. Note \ref{note:injection}), it is of interest to
consider the following algebra.

\begin{definition}
The {\rm higher order peak algebra},  or {\rm generalized peak algebra of order $N$}
is the $\C$-algebra
$\SymN=\Theta_{\zeta}(\Sym)$.
\end{definition}

We will see (Theorem \ref{theorem:complete_basis}) 
that these algebras depend only  on $N$ (and not on
the particular choice of $\zeta$). Note  that
$\SymN$ has obviously a graded algebra structure, inherited
from  $\Sym$
$$
\SymN = \ \bigoplus_{n \geq 0} \ \Sym_n(N)\, ,
$$
where  $\Sym_n(N) := \SymN \cap \Sym_n${\,}.

\begin{note}
{\rm
 For $N = 1$,  $\Sym(1)$ reduces to
 \Sym. Indeed,}
\begin{equation}
\Theta_{q}(S_n) = \ \sum_{i=0}^{n-1}\ (-q)^i \, R_{1^i,n-i}
\end{equation}
{\rm for every $n \geq 1$ (cf. Proposition 5.2 of \cite{NCSF2}),
which immediately implies that $\Theta_{1}(S_n) = \Psi_n$
(cf. Corollary 3.14 of \cite{NCSF1}). 
}
\end{note}

\begin{note}
\textnormal{
For $N = 2$, due to the results of Section \ref{subsec:rel_peak_algebra},
$\Sym(2)$ reduces to $\mathcal{P}$, the inverse image under $\alpha$ of
the usual peak algebra $\Pi$ within $\Sym$.
}
\end{note}

Let us now consider the multiplicative basis $(\pi^I(q))_I$ of
$\Sym$, 
associated with the family of Lie idempotents $(\pi_n(q))_{n \geq
1}$ defined by Formula (\ref{equ:charproppi}) as in Section
\ref{subsec:the_1-q_transformation}. Considering the image of
$(\pi^I(q))_I$ in $\SymN$ by the morphism $\Theta_{\zeta}$, we
obtain that the non-zero elements of the family $(\Theta_{\zeta}(\pi^I(q)))_I$
form a basis of $\SymN$. But according to equation
(\ref{equ:charproppi}), we have that $\Theta_{\zeta}(\pi_n(q)) =
0$ if and only if $n \equiv 0$ $[N]$. Hence the family
\vspace{-0.14cm}
\begin{equation}
\Theta_{\zeta}(\pi^I(q)) = \Theta_{\zeta}(\pi_{i_1}(q))\
\Theta_{\zeta}(\pi_{i_2}(q))\ \cdots\
\Theta_{\zeta}(\pi_{i_r}(q))\, ,
\end{equation}
\noindent indexed by the compositions $I$ of the set $\Ff_n^{(N)}
:= \{\,I = (i_1,\cdots,i_r) \models n \mid \forall\, k \in [1,r],
\, i_k \not\equiv 0$ $[N]\,\}$ forms a basis of $\SymN$.
Summarizing, we obtain the following proposition.

\begin{proposition}
\label{prop_iso}
Let $\zeta$ and $\eta$ be two primitive $N$-th roots of
unity. Then the two $\C$-algebras $\Theta_{\zeta}(\Sym)$ and
$\Theta_{\eta}(\Sym)$ are isomorphic. Moreover their common
dimension is 
$$
dim(\Sym_n(N)) = \#\Ff_n^{(N)} \ .
$$
\end{proposition}
\qed

\paragraph{}
At this stage, it is  interesting to calculate the dimension of
$\Sym_n(N)$.

\begin{proposition}
The Hilbert series of $\SymN$ is
\begin{equation} \label{equ:gen_func}
\sum_{n\ge 0}\,\,dim(\Sym_n(N))\,\,t^n =
\frac{1-t^N}{1-t-t^2- \cdots - t^N}\ .
\end{equation}
\end{proposition}

\Proof The right-hand side of  (\ref{equ:gen_func}) is
the generating function of the generalized Fibonacci numbers
$(f_n)_{n \geq 0}$. It is easily seen to be equal to the generating function
of $\#\Ff_n^{(N)}$, which is
\begin{equation}
\sum_{n\ge 0}\#\Ff_n^{(N)}t^n=\left(1-\sum_{j\not\equiv 0\mod N}t^j\right)^{-1}\,.
\end{equation}
\qed

\begin{note}
{\rm Observe that  $\Sym_n(N) = \Sym_n$ for  $n
\in [1,N\!-\!1]$.}
\end{note}

Though the family $(\Ff_{n}^{(N)})_{n \geq 0}$ was useful for the
calculation of the graded dimension of
$\SymN$, the associated compositions are not  always well-suited
as labellings of the  combinatorial
bases that we want to  construct. 
The following result provides an alternative labelling scheme.

\begin{proposition}
\label{prop:FG} Let us define
$$
\G_{n}^{(N)} := \{\,I = (i_1,\cdots,i_r) \models n \mid \forall\ k
\in [1,r\!-\!1],\, i_k \in [1,N] \ \hbox{and} \  i_r \in [1,
N-1]\,\}.
$$
Then, there exists a bijection between $\G_{n}^{(N)}$ and
$\Ff_{n}^{(N)}$.
\end{proposition}

\Proof Observe that $\G_{n}^{(N)}$ is just
the set of all compositions of $n$ that belong to the rational
language $[1,N]^*[1,N\!-\!1]$. In other words, 
$$
\bigsqcup_{n \geq 1}\,\G_{n}^{(N)} = [1,N]^*[1,N\!-\!1] =
(N^*[1,N\!-\!1])^*N^*[1,N\!-\!1] = (N^*[1,N\!-\!1])^+ \ .
$$
This implies that 
$$
\bigsqcup_{n \geq 1}\,\G_{n}^{(N)} = \{\,N^{i_1}j_1\, \cdots\,
N^{i_r}j_r \mid i_1,i_2,\cdots,i_r \geq 0, j_1, \cdots, j_r \in
[1,N\!-\!1]\,\}\, .
$$
The  map $\varepsilon$ defined by 
$$
N^{i_1}j_1\, N^{i_2}j_2\, \cdots \, N^{i_r}j_r \ \,
\stackrel{\varepsilon}{\longleftrightarrow} \ \,
(N\!\times\!i_1\!+\!j_1,N\!\times\!i_2\!+\!j_2,\cdots,N\!\times\!i_r\!+\!j_r),
$$
for  $i_1,i_2,\cdots,i_r \geq 0$ and $j_1, \cdots, j_r \in
[1,N\!-\!1]$ is then a bijection between
$\G_{n}^{(N)}$ and $\Ff_{n}^{(N)}$. \qed

\begin{note}
{\rm
For $N = 2$ (remembering that
$\Sym(2)$ is the usual peak algebra $\Pi$), $\G_{n}^{(N)}$ is
the set of conjugates of all peak compositions of  $n$, as
introduced in Section \ref{subsec:comp_perm}.}
\end{note}


\subsection{Analogues of complete functions in higher order peak algebra}

If $I = (i_1,i_2,\ldots,i_r)$ and $J = (j_1,j_2,\ldots,j_{s})$ are
two compositions of $n$, let us write 
$I\leq_{p^{(N)}_n} J$ if and only if  $J$ can be
obtained from $I$ by one of the following elementary
rewriting rules:
\begin{center}
$
\begin{array}{rl}
\forall k \in [1,r], & \left\{
\begin{array}{l}
i_k \longmapsto (1, i_k - 1)\,, \wtr{3}{4} \\
i_k \longmapsto (2, i_k - 2)\,, \vtr{3} \\
\hspace{0.5cm}\ldots \vtr{3} \\
i_k \longmapsto (N - 1, i_k - N + 1)\,. \vtr{3}
\end{array}
\right.
\end{array}
$
\end{center}
We denote by $\leq_{P^{(N)}_n}$ the 
transitive closure of $\leq_{p^{(N)}_n}$.

\begin{definition}
For  $n \geq 0$, $P^{(N)}_n$ is the poset of all compositions
of $n$ endowed  with the partial ordering $\leq_{P^{(N)}_n}$.
\end{definition}

\setlength{\unitlength}{0.85pt}

\begin{example}
\textnormal{On compositions of $4$, we have   }
\smallskip
\end{example}

\begin{center}
\begin{tabular}{|rc||rc|}
\hline &$N=2$&& $N=3$ \wtr{2}{4}
\\
\hline &&&
\\
\begin{picture}(30,130)(0,0)
\put(5,65){$P^{(2)}_4 = $}
\end{picture}
&
\begin{picture}(140,130)(0,0)

\put(64,120){(4)} \put(12,80){(1,3)} \put(60,80){(2,2)}
\put(110,80){(3,1)} \put(7,40){(1,1,2)} \put(56,40){(1,2,1)}
\put(107,40){(2,1,1)} \put(53,0){(1,1,1,1)}

\put(67,115){\vector(-3,-2){37}}

\put(70,75){\vector(-3,-2){37}} \put(70,75){\vector(3,-2){37}}

\put(22,75){\vector(0,-1){23}}

\put(120,75){\vector(-3,-2){37}}

\put(71,35){\vector(0,-1){23}} \put(120,35){\vector(-3,-2){37}}
\put(22,35){\vector(3,-2){37}}


\end{picture}

&

\begin{picture}(30,130)(0,0)
\put(5,65){$P^{(3)}_4 = $}
\end{picture}
&

\begin{picture}(140,130)(0,0)
\put(64,120){(4)} \put(12,80){(1,3)} \put(60,80){(2,2)}
\put(110,80){(3,1)} \put(7,40){(1,1,2)} \put(56,40){(1,2,1)}
\put(107,40){(2,1,1)} \put(53,0){(1,1,1,1)}

\put(69,115){\vector(-3,-2){37}} \put(69,115){\vector(0,-1){23}}

\put(70,75){\vector(-3,-2){37}} \put(70,75){\vector(3,-2){37}}

\put(22,75){\vector(0,-1){23}} \put(22,75){\vector(3,-2){37}}

\put(120,75){\vector(-3,-2){37}} \put(120,75){\vector(0,-1){23}}

\put(71,35){\vector(0,-1){23}} \put(120,35){\vector(-3,-2){37}}
\put(22,35){\vector(3,-2){37}}

\end{picture}
\\&&& \\\hline
\multicolumn{2}{|l||}{A single rule: $i_k \longmapsto (1,i_k -
1)$\wtr{2}{4}} &
\multicolumn{2}{l|}{Two rules: $i_k \longmapsto (1,i_k - 1)$, $(2,i_k - 2)$} \\
\hline
\end{tabular}
\end{center}

\bigskip

\medskip

\noindent We can now state one of the main results of this
article.

\begin{theorem} \label{theorem:complete_basis}

For $I$ in $\G_{n}^{(N)}$, let us define the
{\rm complete  peak function of order $N$} $\Sigma_I^{(N)}$ 
by
$$
\Sigma_I^{(N)} := \sum_{J \leq_{P^{(N)}_n}I} R_J \ .
$$
Then, the family $(\Sigma_I^{(N)})_{I \in \G_{n}^{(N)}}$ forms a
basis of $\textnormal{\textbf{Sym}}_n(N)$.
\end{theorem}

\begin{example}
{\rm  }
$$
\begin{array}{l}
\Sigma^{(2)}_{(1,2,1)} = R_{(1,2,1)} + R_{(3,1)} \,, \\
\Sigma^{(3)}_{(1,2,1)} = R_{(1,2,1)} + R_{(3,1)} + R_{(1,3)} + R_{(4)} \,, \\
\\
\Sigma^{(3)}_{(1,1,2)} = \Sigma^{(2)}_{(1,1,2)} = R_{(1,1,2)} +
R_{(2,2)} + R_{(1,3)} + R_{(4)} \,.
\end{array}
$$
\end{example}

\noindent We shall deduce Theorem \ref{theorem:complete_basis}
from  the following two Lemmas.

\begin{lemma} \label{lemma:rnij_elmt_symn}
For all $i \geq 0$ and  $j \in [1, N\!-\!1]$, 
$$
R_{N^ij} \in \SymN \,,
$$
with $N^ij = (\underbrace{N,\ldots\,,N}_{i \ \hbox{\scriptsize
times}},j)$.
\end{lemma}

\begin{note} \label{note:first_elmt}
{\rm  Since $N^ij$ cannot have any predecessor in
$P^{(N)}_n$,  we 
have, for  $i \geq 0$ and  $j \in [1,N\!-\!1]$,}
$$
\Sigma^{(N)}_{N^ij} = R_{N^ij}\,.
$$
\end{note}

\begin{lemma} \label{lemma:product_complete_basic}
For  $I, J \in \G^{(N)}_n$, 
$$
\Sigma^{(N)}_I \times \Sigma^{(N)}_J = \Sigma^{(N)}_{I\cdot J} \, .
$$
\end{lemma}

\vspace{0.5cm} \noindent {\it Proof of Theorem
\ref{theorem:complete_basis} - } Let
$I=(N^{i_1},{j_1},\ldots\,,N^{i_r},{j_r})$ be a generic element of
$\G^{(N)}_n$. According to Lemma
\ref{lemma:product_complete_basic}, we can  write
$$
\Sigma^{(N)}_I = \Sigma^{(N)}_{N^{i_1}j_1} \times \ldots \times
\Sigma^{(N)}_{N^{i_r}j_r}\,\,.
$$
Using now Lemma \ref{lemma:rnij_elmt_symn} and Note
\ref{note:first_elmt}, we deduce from this identity that
$\Sigma^{(N)}_I \in \SymN$. But
the $\Sigma^{(N)}_I$ are linearly independent by construction.
Hence, they form a basis of $\SymN$, since the
cardinality of $\G^{(N)}_n$ is equal to the dimension of
$\Sym_n(N)$, according to Propositions
\ref{prop_iso} and \ref{prop:FG}. \qed

\vspace{0.5cm} \noindent It remains to establish the  Lemmas.

\vspace{0.5cm} \noindent {\it Proof of Lemma
\ref{lemma:rnij_elmt_symn} - }
For $j \in [1,N\!-\!1]$, let
\begin{equation} \label{equ:gen_func_rnij}
\varrho_j(z) \, = \, \sum_{m \geq 0} \, (-1)^{m+1} \, R_{N^mj} \,
z^{mN+j} \,\,.
\end{equation}
Substituting 
$$
R_{N^mj} = \sum_{\scriptsize{\begin{array}{c}
    \alpha_1+\ldots+\alpha_r = m\\
    \alpha_1,\ldots,\alpha_r \geq 0
  \end{array}}} \,
(-1)^{m+1-r} \, S_{\alpha_1N} \, \ldots \, S_{\alpha_{r-1}N} \,
S_{\alpha_{r}N + j}
$$
in (\ref{equ:gen_func_rnij}) yields 
$$
\varrho_j(z) \, = \, \sum_{n \geq 0} \,
\sum_{\scriptsize{\begin{array}{c}
    \alpha_1+\ldots+\alpha_r = n\\
    \alpha_1,\ldots,\alpha_r \geq 0
  \end{array}}} \,
(-1)^r \, S_{\alpha_1N} \, \ldots \, S_{\alpha_{r-1}N} \,
S_{\alpha_rN} \, z^{nN+j} \, \, .
$$
Rearranging the sum, we obtain
$$
\varrho_j(z) \, = \, \sum_{m \geq 0} \, \bigg{(} \,
\sum_{\scriptsize{\begin{array}{c}
    \alpha_1 + \ldots + \alpha_s \geq 0\\
    s \geq 0
  \end{array}}} \,
(-1)^s \, S_{\alpha_1N}\, \ldots \, S_{\alpha_sN} \,
z^{(\alpha_1+\ldots+\alpha_s)N} \bigg{)} \, S_{mN+j} \, z^{mN+j}
\, \, ,
$$
where the two summations are now independent. So, Equation
(\ref{equ:gen_func_rnij}) reduces to 
\begin{equation} \label{equ:}
\varrho_j(z) \, = \, \bigg{(} \,
\sum_{\scriptsize{\begin{array}{c}
    \alpha_1 + \ldots + \alpha_s \geq 0\\
    s \geq 0
  \end{array}}} \,
(-1)^s \, S_{\alpha_1N}\, \ldots \, S_{\alpha_sN} \,
z^{(\alpha_1+\ldots+\alpha_s)N} \bigg{)} \times \bigg{(} \,
\sum_{m \geq 0} \, S_{mN+j} \, z^{mN+j} \bigg{)} \,\, ,
\end{equation}
which can  itself be rewritten as
\begin{equation} \label{equ:rho_ind_prod_inv}
\varrho_j(z) \, = \, \bigg{(} \, \sum_{n \geq 0} \, S_{nN} \,
z^{nN} \bigg{)}^{-1} \times \bigg{(} \, \sum_{m \geq 0} \,
S_{mN+j} \, z^{mN+j} \bigg{)} \, \, .
\end{equation}
It follows then from (\ref{equ:rho_ind_prod_inv}) that
$$
\begin{array}{rl}
\displaystyle{\sum^{N\!-\!1}_{j=1}} \, \varrho_j(z) & = \,
\bigg{(} \, \displaystyle{\sum_{n \geq 0}} \, S_{nN} \, z^{nN}
\bigg{)}^{-1} \times \bigg{(} \, \displaystyle{\sum_{m \geq 0}} \,
S_{mN+1} \, z^{mN+1} \, + \, \ldots \, + \, \displaystyle{\sum_{m
\geq 0}} \, S_{mN+N-1} \, z^{mN+N-1} \bigg{)} \\ & = \, \bigg{(}
\, \displaystyle{\sum_{n \geq 0}} \, S_{nN} \, z^{nN}
\bigg{)}^{-1} \times \bigg{(} \, \displaystyle{\sum_{m \not\equiv
0 \, [N]}} \, S_{m} \, z^{m} \bigg{)} \\ & = \, \bigg{(} \,
\displaystyle{\sum_{n \geq 0}} \, S_{nN} \, z^{nN} \bigg{)}^{-1}
\times \bigg{(} \, \displaystyle{\sum_{m \geq 0}} \, S_{m} \,
z^{m} - \displaystyle{\sum_{n \geq 0}} \, S_{nN} \, z^{nN}
\bigg{)} \\ & = \, \bigg{(} \, \displaystyle{\sum_{n \geq 0}} \,
S_{nN} \, z^{nN} \bigg{)}^{-1} \times \bigg{(} \,
\displaystyle{\sum_{m \geq 0}} \, S_{m} \, z^{m} \bigg{)} \, - \,
1 \,,
\end{array}
$$
so that
\begin{equation} \label{}
1 \, + \, \sum^{N\!-\!1}_{j=1} \, \varrho_j(z) \, = \, \bigg{(} \,
\sum_{n \geq 0} \, S_{nN} \, z^{nN} \bigg{)}^{-1} \times \bigg{(}
\, \sum_{m \geq 0} \, S_{m} \, z^{m} \bigg{)} \,.
\end{equation}
%

We want to prove that this series in $z$ has its coefficients in $\SymN$.
Equivalently, we can prove this for its inverse
\begin{eqnarray}
\left(1 \, + \, \sum^{N\!-\!1}_{j=1} \, \varrho_j(z)\right)^{-1}
&= \displaystyle \lambda_{-z}(A)\sum_{m\ge 0}S_{mN}(A)z^{mN}\nonumber\\
&= \displaystyle \sum_{n\ge 0}z^n \sum_{Ni+j=n}(-1)^j\Lambda_j(A)S_{Ni}(A)\,.
\end{eqnarray}
The coefficient $C_n$ of $z^n$ in this expression can be rewritten as
\begin{equation}
\sum_{Ni+j=n}(-1)^j S_j(-A)S_{Ni}(qA)|_{q=\zeta} =E_0[S_n((q-1)A)]
\end{equation}
where for a polynomial $f(q)$ we set
\begin{equation}
f(q)= f_0(q^N)+q f_1(q^N)+q^2f_2(q^N)+\cdots+q^{N-1}f_{N-1}(q^N)
\end{equation}
and 
\begin{equation}
f_j = E_j[f] \,.
\end{equation}
Now, the inversion formula for the discrete Fourier transform on $N$
points shows that the polynomials $q^jf_j(q^N)$ are linear combinations
with complex coefficients of $f(q)$, $f(\zeta q)$, $f(\zeta^2q)$,...,$f(\zeta^{N-1}q)$.
In particular, $E_0[S_n((q-1)A)]$ is a linear combination of
$$S_n((q-1)A),\ S_n(\zeta q-1)A,\ldots, S_n(\zeta^{N-1}q-1)A)\,,$$
and $C_n$, which is its specialization at $q=\zeta$, is therefore
in the subspace spanned by
$$
S_n((\zeta-1)A),\ S_n((\zeta^2-1)A),\ldots,\     S_n((\zeta^N-1)A)=0\,
$$
which are all in $\SymN$, thanks to the identity
$\zeta^k-1=(1-\zeta)(-1-\zeta-\cdots -\zeta^{k-1})$.
\qed

\vspace{0.5cm} \noindent {\it Proof of Lemma
\ref{lemma:product_complete_basic} - } Observe first that for $I,
J \in \G^{(N)}_n$, the concatenation $I\cdot J$ is also in
$\G^{(N)}_n$. Consider now the product
\begin{equation} \label{equ:product_r}
\bigg{(} \, \sum_{K \leq_{P^{(N)}_n} I} \, R_K \bigg{)} \times
\bigg{(} \, \sum_{L \leq_{P^{(N)}_n} J} \, R_L \bigg{)} \,\,.
\end{equation}
According to the product formula for ribbons 
(Proposition 3.13 of \cite{NCSF1}), the product
in Formula (\ref{equ:product_r}) is seen to be the sum of all
compositions of $n$ which can be refined into the concatenation
$I.J$, whence the Lemma.
\qed

\begin{note}
{\rm One should observe that the converse of Lemma
\ref{lemma:rnij_elmt_symn} is false:  there exists
$R_I \in \SymN$ such that the associated composition $I$ is
different from $N^ij$ for all $i \geq 0$ and $j \in
[1,N\!-\!1]$. For example,}
$$
R_{(2,1,1)} \, = \, \bigg{(} \, \Sigma^{(3)}_{(2,1,1)} -
\Sigma^{(3)}_{(3,1)} - \Sigma^{(3)}_{(2,2)} \bigg{)} \, \in \SymN
\,,
$$
{\rm and $(2,1,1)$ is obviously not a composition of the type
$(3^i,1)$ or $(3^i,2)$.}
\end{note}

\vspace{0.5cm} \noindent For  $k \not\equiv 0 \, [N]$, let us
now set
$$
T_k = R_{N^ij} = \Sigma_{N^ij}
$$
where $(i,j)$ is the unique pair of $\N \times [1, N\!-\!1]$ such
that $k = N \times i + j$. More generally, for 
$K = (k_1,\ldots,k_r) \in \Ff^{(N)}_n$, we set
$$
T_K = T_{k_1} \times \ldots \times T_{k_r} \,\,.
$$
An alternative definition of $T_K$ is 
$$
\Sigma_I = T_{\varepsilon(I)} \,\, ,
$$
where $\varepsilon$ is the bijection introduced in the proof of
Proposition \ref{prop:FG}, which allows us to
rephrase Theorem \ref{theorem:complete_basis} in the form

\begin{corollary}
The family $(T_K)_{K \in \Ff^{(N)}_n}$ forms a basis of
$\Sym^{(N)}_n$.
\end{corollary}


\subsection{Generalized peak ribbons}

We are now in a position to introduce analogues of ribbons
in the higher order peak algebras.

\begin{definition} \label{def:rho}
Let $I$ be a composition of $\G^{(N)}_n$. Then the {\rm peak ribbon of order $N$}
labelled by $I$ is defined as
\begin{equation} \label{equ:decomp_rho_on_sigma}
\rho^{(N)}_I \, = \, \sum_{\scriptsize{\begin{array}{c}
      J \leq_{P^{(N)}_n} I \,,\\
      J \in \G^{(N)}_n
  \end{array}}} \,
(-1)^{l(I) - l(J)} \, \Sigma^{(N)}_J \,.
\end{equation}
\end{definition}

\noindent 
The following
Proposition results  immediately from Theorem
\ref{theorem:complete_basis}.

\begin{proposition}
  $(\rho^{(N)}_I)_{I \in \G^{(N)}_n}$ is a basis of $\SymN$.
\end{proposition}

\noindent It follows that each
$\Sigma^{(N)}_I$ has a decomposition on the 
ge\-ne\-ra\-li\-zed peak ribbons. This decomposition is given by
the M\"obius-like inverse of Formula
(\ref{equ:decomp_rho_on_sigma}), justifying the claim that the families
$(\rho^{(N)}_I)_{I \in \G^{(N)}_n}$ and $(\Sigma^{(N)}_I)_{I \in
\G^{(N)}_n}$ are in a similar relation as  $(R_I)_{I \models n}$ and $(S_I)_{I
\models n}$.

\begin{proposition} \label{prop:decRS}
For all $I \in \G^{(N)}_n$, 
\begin{equation} \label{equ:decomp_sigma_on_rho}
  \Sigma^{(N)}_I
  \, = \,
\sum_{\scriptsize{\begin{array}{c}
      J \leq_{P^{(N)}_n} I \\
      J \in \G^{(N)}_n
  \end{array}}} \,
\rho^{(N)}_J
\end{equation}
\end{proposition}

\Proof See Note \ref{proof:RS}. \qed

\begin{note}
{\rm We can introduce others analogues of the usual 
ribbon functions, such as, for example, the 
following two families, which are indexed by compositions $I \in
\G^{(N)}_n$:}
$$
\rho^{(N)}_I(t) \,= \,
\displaystyle{\sum_{\scriptsize{\begin{array}{c}
    J \leq_{P^{(N)}_n} I\\
    J \in \G^{(N)}_n
\end{array}}}} \,
t^{l(I) - l(J)} \, \Sigma_J \, \quad \ \hbox{{\rm or}} \quad \
\rho'^{(N)}_I(t) \,= \,
\displaystyle{\sum_{\scriptsize{\begin{array}{c}
      J \leq_{P^{(N)}_n} I\\
      J \in \G^{(N)}_n
    \end{array}}}} \,
    t^{l(I) + l(J)} \,
    \Sigma_J
    \, ,
\wtr{3}{6}
$$
{\rm and which also are clearly bases of $\SymN$ for  $t \in
\C - \{0\}$. We have of course
$$
\rho^{(N)}_I = \rho^{(N)}_I(-1) = \rho'^{(N)}_I(-1) \,,
$$
for  $I \in \G^{(N)}_n$. Note also that $\rho^{(N)}_I(t)$ and
$\rho'^{(N)}_I(t)$ are connected by}
$$
\rho'^{(N)}_I(t) \, = \, t^{2l(I)} \, \rho^{(N)}_I(\frac{1}{t}) \
.
$$
{\rm Due to this relation, we will only make use of
$(\rho^{(N)}_I(t))_{I \in \G^{(N)}_n}$ in the sequel.}
\end{note}

\begin{example}
{\rm We now give explicitly some elements of the family
$(\rho^{(N)}_I)_{I \in \G^{(N)}_n}$: }
$$
\rho^{(2)}_{(1,1,1)} \, = \, \Sigma^{(2)}_{(1,1,1)} -
\Sigma^{(2)}_{(2,1)}\, , \quad \left\{\
\begin{array}{l}
\rho^{(3)}_{(1,1,1)} \, = \, \Sigma^{(3)}_{(1,1,1)} -
\Sigma^{(3)}_{(1,2)} - \Sigma^{(3)}_{(2,1)} \, = \, R_{(1,1,1)} -
R_{(3)} \,,
\\
\rho^{(3)}_{(1,2)} \, = \, \Sigma^{(3)}_{(1,2)} \, = \, R_{(1,2)}
+ R_{(3)} \,,
\\
\rho^{(3)}_{(2,1)} \, = \, \Sigma^{(3)}_{(2,1)} \, = \, R_{(2,1)}
+ R_{(3)} \,.
\end{array}
\right.
$$
\end{example}

\noindent Unfortunately the internal product of two elements of
the family $(\rho^{(N)}_I(t))_{I \in \G^{(N)}_n}$ does not always
decompose with non-negative coefficients on this family, contrary to
the case of the usual ribbon basis of $\Sym$ (cf. section 5 of
\cite{NCSF1}). This shows that the family defined by Definition
\ref{def:rho} is probably not a perfect analogue of the usual non
commutative ribbons.

\begin{example}{\rm }
  $$
  \left\{
  \begin{array}{l}
    \rho^{(3)}_{(1,1,1)} \ast \rho^{(3)}_{(1,1,1)}
    \, = \,
    (-2) \times \rho^{(3)}_{(1,1,1)} \,,
    \\
    \rho^{(3)}_{(1,2)} \ast \rho^{(3)}_{(1,1,1)}
    \, = \,
    \rho^{(3)}_{(1,1,1)} + \rho^{(3)}_{(2,1)} - \rho^{(3)}_{(1,2)} \,.
  \end{array}
  \right.
  $$
\end{example}


\section{Decompositions on peak bases}

\subsection{A projector}

Let us define a linear map $\pi_N$ from $\Sym$ onto $\SymN$
by 
\begin{equation}
\pi_N(S^I) \, = \, \left\{
 \begin{array}{ll}
   \Sigma_I
   &
   $if $ I \in \G^{(N)}_n\,,
   \\ \\
   0
   &
   $otherwise$.
 \end{array}
 \right.
\end{equation}
\begin{note} \label{note:rhoi_ri}
{\rm Note that $\pi_N(R_I) = \rho^{(N)}_I$ for every $I \in
\G^{(N)}_n$.}
\end{note}

Let $T^{(N)}$ be the left-ideal of $\Sym$ generated by the $S_j$ such that
$j\not\equiv 0 \mod N$.

\begin{lemma} \label{lemma:tn_include_symn}
$\SymN$ is a subalgebra of $T^{(N)}$.
\end{lemma}

\Proof 
The 
generators $R_{N^ij}$ of the algebra $\SymN$ are  in
$T^{(N)}$, since
$$
R_{N^ij} = \sum_{\scriptsize{\begin{array}{c}
      \alpha_1 + ...\ldots + \alpha_r = i\\
      \alpha_1, \ldots , \alpha_r \geq 0
\end{array}}} \,
(-1)^{i+1-r} \, S^{(\alpha_1N, \, \ldots \,,\, \alpha_{r-1}N,\,
\alpha_{r}N + j)} \,.
$$
\qed

\begin{proposition} \label{prop:pin_alg_morphism}
$\pi_N|_{\SymN}$ is an algebra morphism.
\end{proposition}

\Proof We will in fact show the slightly stronger property
\begin{equation} \label{equ:morph_algebra_larger_stat}
\pi_N(F \times G) \, = \, \pi_N(F) \times \pi_N(G) \,,
\end{equation}
for all $F \in T^{(N)}$ and $G \in \Sym$. 
We  establish 
(\ref{equ:morph_algebra_larger_stat}) by proving that
it holds for products of complete
functions, that is
\begin{equation} \label{equ:morph_algebra_right_stat}
\pi_N(S^I \times S^J) \, = \, \pi_N(S^I) \times \pi_N(S^J) \,,
\end{equation}
for  $S^I \in T^{(N)}_n$ and all
$S^J$. Two cases are to be considered.

\smallskip
Suppose first that $I\cdot J \not\in \G^{(N)}_n$. Then one has either
$I \not\in \G^{(N)}_n$, or $J \not\in \G^{(N)}_n$ by definition of
$\G^{(N)}_n$. So, both sides of 
(\ref{equ:morph_algebra_right_stat}) are  $0$
according to the definition of $\pi_N$.

\smallskip
Assume now that $I\cdot J \in \G^{(N)}_n$. Then, $I$ and
$J$ are also elements of $\G^{(N)}_n$, since we have assumed that
the last part of  $I$ is not a multiple of $N$.  Then, by definition of
$\pi_N$, 
\begin{equation} \label{equ:pin_sij_eq_sigmanij}
\pi_N(S^{I\cdot J}) \, = \, \Sigma^{(N)}_{I\cdot J} \,.
\end{equation}
But from Lemma \ref{lemma:product_complete_basic}, we know that
that
\begin{equation} \label{equ:sigmanij_eq_sigmani_sigmanj}
\Sigma^{(N)}_{I.J} \, = \, \Sigma^{(N)}_I \times \Sigma^{(N)}_J
\,.
\end{equation}
So it results from (\ref{equ:pin_sij_eq_sigmanij}) and
(\ref{equ:sigmanij_eq_sigmani_sigmanj}) that 
$$
\pi_N(S^{I\cdot J}) \, = \, \pi_N(S^I) \times \pi_N(S^J) \,,
$$
since $\Sigma^{(N)}_I$ and $\Sigma^{(N)}_J$ are respectively the
images of $S^I$ and $S^J$ under $\pi_N$. Hence, 
(\ref{equ:morph_algebra_right_stat}) also holds in this case.
\qed

\vspace{0.5cm} \noindent We can now give an important property of $\pi_N$.

\begin{proposition}
$\pi_N$ is a projector from $\Sym$ onto $\SymN$.
\end{proposition}

\Proof We have to prove that
\begin{equation} \label{equ:projector_formula}
\pi_N(\Sigma_I) \, = \, \Sigma_I \,,
\end{equation}
for $I \in \G^{(N)}_n$. Using Proposition
\ref{prop:pin_alg_morphism}, which gives that $\pi_N$ is
multiplicative on $\SymN$, and Lemma
\ref{lemma:product_complete_basic}, we only have to prove that
$$
\pi_N(R_{N^ij}) = R_{N^ij} \,,
$$
for  $i \geq 0$ and $j \in [1,N\!-\!1]$. But we have already seen (Note \ref{note:rhoi_ri})
that 
$\pi_N(R_I) = \rho^{(N)}_I$ for every $I \in
\G^{(N)}_n$, and for $I=(N^ij)$, $\rho^{(N)}_I=\Sigma_I=R_I$.
\qed

\begin{note} \label{proof:RS}
{\rm Remark that for  $I \in \G^{(N)}_n$,}
$$
\Sigma^{(N)}_I \, = \, \pi_N(S^I) \, = \, \sum_{J \leq I} \,
\pi_N(R_J) \, = \, \sum_{\scriptsize{\begin{array}{c}
      J \leq_{P^{(N)}_n} I\\
      J \in \G^{(N)}_n
  \end{array}}} \,
\rho_J \,,
$$
{\rm which gives an alternative proof of Proposition
\ref{prop:decRS}.}
\end{note}


\subsection{Some interesting decompositions}

We record in this section a number of remarkable decompositions
with respect to the bases of  
complete and ribbon peak
functions of order $N$.

\begin{proposition} \label{prop:decomp_thetasi_sigmaj}
Let $\zeta$ be a primitive $N$-th root of unity. Then, for a 
composition $I$ of $n$, 
\begin{equation} \label{equ:thetazetasi_eq_sum_sigmaj}
\vartheta_{\zeta}(S^I) \, = \, \sum_{\scriptsize{\begin{array}{c}
 J=(j_1,...,j_s) \leq I \\
 J \in \G^{(N)}_n
\end{array}}} \,
(-1)^{l(I) - l(J)} \, \prod_{k \in \H(I,J)} \, \bigg{(}
\bigg(\frac{1}{\zeta}\bigg)^{j_k} - 1 \bigg{)} \, \Sigma^{(N)}_J
\,,
\end{equation}
where 
$$
\H(I,J) \, := \, \{ \, l \in [1,s], | \exists k \in [1,r],
j_1+...+j_l=i_k \, \} \,.
$$
\end{proposition}

\Proof Formula (\ref{equ:thetazetasi_eq_sum_sigmaj}) results by
applying $\pi_N$ to the decomposition of $S^I((1-q)A)$ on products
of complete functions,
specialized at $q = \zeta$, 
which is given by Formula
(105) of Proposition 5.30 of \cite{NCSF2}. \qed

\begin{proposition} \label{prop:decomp_thetari_sigmaj}
Let $\zeta$ be a primitive $N$-th root of unity. Then, for a 
composition $I$ of $n$, 
\begin{equation} \label{equ:thetazetari_eq_sigmaj}
\vartheta_{\zeta}(R_I) \, = \, \sum_{\scriptsize{\begin{array}{c}
 J=(j_1,...,j_s) \leq I \\
 J \in \G^{(N)}_n
\end{array}}} \,
(-1)^{l(I) - l(J)} \, q^{\alpha(I,J)} \, (1-q^{j_s}) \,
\Sigma^{(N)}_J \,,
\end{equation}
where 
$$
\alpha(I,J) \, =\, \sum^{l(J)-1}_{k=1} \, j_k \times
\delta_{j_1+...+j_k \not\in D(I)} \,.
$$
\end{proposition}

\Proof As for Proposition \ref{prop:decomp_thetasi_sigmaj},
Formula (\ref{equ:thetazetari_eq_sigmaj}) is obtained by
applying $\pi_N$ on the decomposition of $R_I((1-\zeta)A)$ on the
basis $S^J$, as given in \cite{NCSF2} . \qed

\begin{note}
{\rm Note that Propositions \ref{prop:decomp_thetasi_sigmaj} or
\ref{prop:decomp_thetari_sigmaj} give another proof of Theorem
\ref{theorem:complete_basis}, since Formulas
(\ref{equ:thetazetasi_eq_sum_sigmaj}) or
(\ref{equ:thetazetari_eq_sigmaj}) show that $(\Sigma^{(N)}_I)_{I
\in \G^{(N)}_n}$ is a generating family of $\SymN$. Since this
family is obviously linearly free by construction, it follows that
it is a basis of $\SymN$}
\end{note}

\noindent In order to state the next Proposition, we need  a
definition. For each pair of  compositions
$I=(i_1,...,i_r)$ and $J=(j_1,...,j_s)$  of the same weight, let
us  introduce the sequence of compositions $H(I,J) = (H_1,
... , H_r)$ of length $l(I)$ which is uniquely determined by the
the conditions $|H_k| = i_k$ for  $k \in [1, l(I)]$, and 
$$
H_1 \bullet \, ... \, \bullet H_r \, = \, J \,,
$$
where $\bullet$ denotes either  
the concatenation of compositions, or the operation
$\rhd$, 
defined by 
$$
H \rhd K \, = \, (k_1,..., k_{r-1}, k_r + l_1, l_2, ..., l_s) \,.
$$

\setlength{\unitlength}{1pt}

\begin{example}
{\rm Let $I=(3,2,1,4)$ and $J=(2,5,2,1)$. Then, we have
$$
H(I,J) = ((2,1),(2),(1),(1,2,1)) \,.
$$
%
%
}
\end{example}

\noindent We can now  give the following definition.

\begin{definition}
For any two compositions $I$ and $J$ of the same weight, we set
\begin{equation}
h(I,J) \, = \, \left\{
\begin{array}{ll}
  - \infty &
  {\rm if\: there\: exists}\: k \in [1,l(I)],
  \alpha \geq 0,
  \beta \geq 1,
  H_k \not= (1^{\alpha},\beta) \,,
  \\ \\
  \displaystyle{\sum^{l(I)}_{i=1}} \, \alpha_i &
  {\rm if\: for\: every}\: k \in [1,l(I)],
  \alpha_k \geq 0,
  \beta_k \geq 1,
  H_k = (1^{\alpha_k},\beta_k) \,.
\end{array}
\right.
\end{equation}
\end{definition}

\begin{proposition}
Let $\zeta$ be a primitive $N$-th root of unity. Then, for a
composition $I$ of $n$, 
\begin{equation} \label{equ:decomp_s_rho}
\vartheta_{\zeta}(S^I) \, = \, (1-\zeta)^{l(I)} \,
\sum_{\scriptsize{\begin{array}{c}
    J \in \G^{(N)}_n \\
    h(I,J) \not= -\infty
\end{array}}} \,
(-\zeta)^{h(I,J)} \, \rho^{(N)}_J \,
\end{equation}
\end{proposition}

\Proof The decomposition of $\vartheta_{q}(S^I)$ on the ribbon basis 
follows from Formula (67) of
Proposition 5.2 of \cite{NCSF2}. Formula (\ref{equ:decomp_s_rho})
comes by applying $\pi_N$ to this decomposition, specialized
at $q=\zeta$. \qed

\vspace{0.5cm} \noindent Note  that an arbitrary composition
$I$ can always  be uniquely written as
$$
I \, = \, H_1 \cdot  H_2 \cdot \cdots \cdot H_{hl(I)} \,,
$$
where $H_k = (1^{\alpha}, \beta)$,  $k \in [1, hl(I)]$.
We denote by $H_I$ the composition
$$
H_I \, = \, (|H_1|, |H_2|, \, ... \, , |H_{hl(I)}|) \,.
$$

\begin{example}
{\rm For $I = (1,3,1,4,2)$, we have 
$$
I =  H_1\cdot  H_2 \cdot H_3 \,,
$$
with $H_1 = (1,3)$, $H_2 = (1,4)$, $H_3 = (2)$.}
%
%
%
%
%
%
%
%
{\rm Hence $hl(I) = 3$ and $D(H_I) = \{\, 4,9 \,\}$, in this
case.}
\end{example}

\begin{definition}
Let $I ,J$ be two compositions of $n$. We set
\begin{equation}
b(I,J) \, = \, \left\{
\begin{array}{cl}
  |(1 + (D(I) - D(J))) \cup (D(J) - D(I))| &
  {\rm if}\: D(H_J) \subset  S(I)\,,
  \\ \\
  - \infty &
  {\rm otherwise}\,,
\end{array}
\right.
\end{equation}
where $S(I) \, := \, ((1 + D(I)) - D(I)) \cup (D(I) - (1+D(I)))$.
\end{definition}

\begin{proposition}
Let $\zeta$ be a primitive $N$-th root of unity. Then, for a
composition $I$ of $n$, 
\begin{equation} \label{equ:decomp_r_rho}
\vartheta_{\zeta}(R_I) \, = \, \sum_{\scriptsize{\begin{array}{c}
    J \in \G^{(N)}_n \\
    b(I,J) \not= -\infty
\end{array}}} \,
(1-\zeta)^{hl(J)} \, (-\zeta)^{b(I,J)} \, \rho^{(N)}_J \,.
\end{equation}
\end{proposition}

\Proof This follows  by applying $\pi_N$ to the ribbon
decomposition of $\vartheta_{q}(R_I)$, specialized at $q=\zeta$,
as given in Formula (121) of Proposition 5.41 of \cite{NCSF2}. \qed


\subsection{Higher order ``noncommutative tangent numbers''}

We will finally present in this last section, a number of formulas
generalizing the known relations involving the so-called 
{\em noncommutative tangent
numbers}, introduced in \cite{NCSF1}. The formulas of
\cite{NCSF1} are obtained for $N = 1$, and those of \cite{BHT}
for $N = 2$.

\begin{proposition} \label{prop:tangent_geom}
Let 
$$
\t \, := \, \sum_{i \geq 0} \, (-1)^{i+1} \, \bigg( \,
\sum_{j=1}^{N-1} \, \Sigma^{(N)}_{N^ij} \bigg)\,,
$$
Then,
\begin{equation} \label{equ:power_serie_rho1n}
(1-\t)^{-1} \, = \, \sum_{n\ge 0} \, (-1)^n \,
\rho^{(N)}_{1^n}\,,
\end{equation}
\end{proposition}

\Proof We know that
\begin{equation} \label{equ:1-t-1_first_sum}
(1-t)^{-1} \, = \, \sum_{\scriptsize{\begin{array}{c}
      i_1,...,i_r \geq 0 \\
      r \geq 0
\end{array}}} \,
(-1)^{i_1+...+i_r+r} \, \sum_{j_1,...,j_r \in [1,N\!-\!1]} \,
\Sigma^{(N)}_{N^{i_1}j_1 \, ... \, N^{i_r}j_r} \,,
\end{equation}
where  $\Sigma^{(N)}_{\emptyset} = 1$. But for a composition
$K=(N^{i_1}j_1, ... , N^{i_r}j_r)$, one has
$$
l(K) = i_1 + ... + i_r +r \,,
$$
so that (\ref{equ:1-t-1_first_sum}) can be rewritten as
$$
\begin{array}{rl}
(1-\t)^{-1} & \, = \, 1 +
\displaystyle{\sum_{\scriptsize{\begin{array}{c}
      K \in \G^{(N)}_n \\
      n \geq 1
\end{array}}}} \,
(-1)^{l(K)} \, \Sigma^{(N)}_K \,
\\ &
\, = \, 1 + \displaystyle{\sum_{n\ge 1}} \, (-1)^n \,
\bigg( \, \displaystyle{\sum_{K \in \G^{(N)}_n}} \, (-1)^{l(K) -
|K|} \, \Sigma^{(N)}_K \, \bigg) \,,
\\ &
\, = \, 1 + \displaystyle{\sum_{n\ge 1}} \, (-1)^n \,
\rho^{(N)}_{1^n}\,,
\end{array}
$$
which is the required expression ($\rho^{(N)}_{1^0} = 1$).
\qed

\begin{note}
{\rm  Equation (\ref{equ:power_serie_rho1n})  can be regarded as a kind of
generalization of the trivial relation $\sigma(1) \, \lambda(-1) = 1$,
where $\sigma$ and $\lambda$ are respectively the generating
functions of  complete and elementary non commutative symmetric
functions. More generally, we can set
$$
\left\{
\begin{array}{l}
\sigma_N(t) = 1 + \displaystyle{\sum^{N-1}_{i=1}} \, t^i \,
\Sigma^{(N)}_{i} + \sum_{i \geq 1} \, (-1)^i \, \bigg( \,
\displaystyle{\sum^{N-1}_{j=1}} \, t^{iN+j} \, \Sigma^{(N)}_{N^ij}
\bigg) \,,
\\ \\
\lambda_N(t) = \displaystyle{\sum_{n \geq 0}} \, (-t)^n \,
\rho^{(N)}_{1^n} \,,
\end{array}
\right.
$$
and one can check, following the lines of the proof of Proposition
\ref{prop:tangent_geom}, 
that $\sigma_N(t) \, \lambda_N(-t) = 1$,  another
generalization of the same relation. }
\end{note}

\begin{proposition}
Let $\zeta$ be a primitive $N$-th root of unity and let
$\t_{\zeta}$ be defined by 
$$
t_{\zeta} \, := \, \sum_{i \geq 0} \, \sum_{j=1}^{N-1} \,
\zeta^{j-i-1} \, \Sigma^{(N)}_{N^ij} \,.
$$
Then,
\begin{equation} \label{equ:power_serie_rho1nzeta}
(1-\t_{\zeta})^{-1} \, = \, \sum_{n\ge 0} \, (-1)^n \,
\rho^{(N)}_{1^n}(\zeta)\,.
\end{equation}
\end{proposition}

\begin{note}
{\rm Note that for $N = 2$, we get $\t_{-1} = TH$,
in the notation of \cite{NCSF1}.
}
\end{note}

\Proof According to Lemma \ref{lemma:product_complete_basic}, 
\begin{equation} \label{equ:1-tz-1_sum_sigma}
(1-\t_{\zeta})^{-1} \, = \, \sum_{\scriptsize{\begin{array}{c}
      i_1,...,i_r \geq 0 \\
      j_1,...,j_r \in [1,N\!-\!1] \\
      r \geq 0
\end{array}}} \,
\zeta^{j_1+...+j_r - (i_1+...+i_r+1)} \,
\Sigma^{(N)}_{N^{i_1}j_1...N^{i_r}j_r} \,.
\end{equation}
But for a composition $K=(N^{i_1}j_1, ... , N^{i_r}j_r)$,
$$
|K| - l(K) \, = \, \underbrace{N\times(i_1 + ... + i_r) + j_1 +
... +j_r}_{|K|} - \underbrace{(r + i_1 + ... +i_r)}_{l(K)}\,. 
$$
So, (\ref{equ:1-tz-1_sum_sigma}) becomes
$$
(1-\t_{\zeta})^{-1} \, = \, \sum_{\scriptsize{\begin{array}{c}
      K \in  \G^{(N)}_n \\
      n \geq 0
      \end{array}}} \,
\zeta^{|K| - l(K)} \, \Sigma^{(N)}_K \,,
$$
whence the proposition. \qed

\bigskip
{\small
\noindent{\it Acknowledgements.-}
This project has been partially supported by EC's IHRP Programme, grant
HPRN-CT-2001-00272, ``Algebraic Combinatorics in Europe".
The  authors are also grateful to Maxime
Rey for his help with the redaction of this paper.

}


\end{document}